# Henri Lebesgue and the End of Classical Theories on Calculus

Xiaoping Ding

*Beijing institute of pharmacology of 21st century*

**Abstract** In this paper a novel calculus system has been established based on the concept of 'werden'. The basis of logic self-contraction of the theories on current calculus was shown. Mistakes and defects in the structure and meaning of the theories on current calculus were exposed. A new quantity-figure model as the premise of mathematics has been formed after the correction of definition of the real number and the point. Basic concepts such as the derivative, the differential, the primitive function and the integral have been redefined and the theories on calculus have been reestablished. By historical verification of theories on calculus, it is demonstrated that the Newton-Leibniz theories on calculus have returned in the form of the novel theories established in this paper.

**Keywords:** Theories on calculus, differentials, derivatives, integrals, quantity-figure model

**Introduction**

The combination of theories and methods on calculus is called calculus. The calculus, in the period from Isaac Newton, Gottfried Leibniz, Leonhard Euler to Joseph Lagrange, is essentially different from the one in the period from Augustin Cauchy to Henri Lebesgue. The sign for the division is the publication of the monograph 'cours d'analyse' by Cauchy in 1821. In this paper the period from 1667 to 1821 is called the first period of history of calculus, when Newton and Leibniz established initial theories and methods on calculus; the period after 1821 is called the second period of history of calculus, when Cauchy denied the common ideas of Newton and Leibniz and established new theories on calculus (i.e., current theories on calculus) using the concept of the limit in the form of the expression invented by Leibniz. The new theories on calculus and the traditional methods on calculus are called mathematical analysis or calculus sometimes.

The basic ideas of theories on calculus in the first period (i.e., Newton-Leibniz theories on calculus) are correct although they are not consistent. Cauchy agreed with formulas ' $dx^e = ex^{e-1}dx$, $\int x^e dx = \dfrac{x^{e+1}}{e+1}$ and $\int_a^b y dx = Z(b) - Z(a)$ ' [1], given by Leibniz, but he did not understand Leibniz's explanation of the differentials $dx$ and $dy$. Leibniz said, 'A differential is like the contact angle of Euclid, which is smaller than any given quantity, but not equal to zero.' He also said, 'We consider an infinitesimal quantity as a relative zero, not a simple zero or an absolute zero.' It should be noted that a worldwide mistake exists that the differential is considered as an arbitrarily small quantity. Indeed, Leibniz indicated that the differential is not zero, or a finite quantity, not to mention infinity, but 'a relative zero smaller than any given quantity'. It is a new type of quantity since the concept of modern numbers hasn't appeared in Leibniz's time; while the concept of 'the relative zero' existed in the concepts of all the numbers, including non-standard numbers of Robinson, from the time of Cantor-Dedekind to 2010. However, Cauchy introduced the idea of the limit to calculus system because he wasn't able to understand the ideas of Leibniz.

Cauchy defined the differential as a finite quantity, which had changed Leibniz's idea that the differential was defined as 'a qualitative zero', 'smaller than any given quantity'. It seems that a derivative has been defined perfectly and the Newton-Leibniz equation has been proved by making use of the concept of limit, but there exist radical problems of equating $dx$ and $\Delta x$ when it comes to defining the differential, because the differential (expressed starting with the symbol 'd') has been defined as the linear main-part of a change (expressed starting with the symbol '$\Delta$') resulting in that $dx$ is not equal to $\Delta x$. However, he continued to establish his 'theories on calculus' producing more associated mistakes.

It is impossible to find out whether the mathematicians, from Bernhard Riemann, Karl Weierstrass to Henri Lebesgue, had discovered the deadly mistakes in the theories on calculus established by Cauchy due to a lack of historical materials. But we are sure that if they revealed the deadly mistakes, what they had done was to continue the mistakes.

The birth of Lebesgue integral led to the establishment of Real Variable Function (i.e., Modern Analysis). Then it has been announced that calculus (mathematical analysis) is a 'rigorous and integrated system' in the mathematical field.

In fact no new concepts on the traditional derivatives, primitive functions (indefinite integrals) or differentials, except the definite integral, were established in the theories in Modern Analysis or Real Variable Function based on the ideas of Lebesgue. Even the Lebesgue integral can be transformed to the Riemann integral in a continuous interval. Therefore, Modern Analysis does not correct the radical mistakes in the current theories on calculus, except the establishment of the Riemann integral and some unimportant new concepts (e.g., variation) as well as the intensification of signifying.

In this paper the defects of the current theories on calculus were shown. The definition of the real number and the point were corrected and a new quantity-figure model has been suggested based on the studies of previous researchers. Basic concepts such as the derivative, the differential, the primitive function and the integral have been redefined and theories on calculus have been reestablished according to the novel definition of differentials based on the concept of 'werden'. The historical verification of the novel theories was carried out by analyzing the history of calculus. The related problems on the establishment of the novel mathematical model were discussed from the viewpoint of mathematical philosophy.

## 1 Mistakes and defects in the current theories on calculus

In the theories on current calculus (i.e., Cauchy-Lebesgue system), the linear main-part (i.e., $f(x_0)\Delta x$) of equation (1.1) is defined as the differential.

$$\Delta y = f(x_0)\Delta x + o(\Delta x) \qquad (1.1)$$

Where $\Delta x$ isn't infinitesimal.

In order to obtain the formula $dy = f(x_0)dx$, two assumptions have been given. One is that the differential $dx$ is considered to be equal to $dy$ based on the function $y = x$ without logic [2,3]; the other is the definition of $dx = \Delta x$ without logic [4-6], either. Therefore, $dy = f(x_0)\Delta x$ is transformed into $dy = f(x_0)dx$ according to one of the two above assumptions.

Either consideration or definition without logic violates the principles of science. Actually, the consideration that the differential $dx$ is equal to the differential of a function $y = x$ leads to not only the required formula $dy = f(x_0)\Delta x = f(x_0)dx$ according to $dx = dy = x' \cdot \Delta x = \Delta x$, but also the ridiculous one $dy = f(x_0)\Delta x = f(x_0)dy$. Regarding a function of two variables $z = f(x, y)$, 'the differential of the independent variable x is considered as the differential of the function $z = x$, thus the formula $dx = dz = \frac{\partial x}{\partial x}\Delta x + \frac{\partial x}{\partial y}\Delta y = 1 \cdot \Delta x + 0 \cdot \Delta y = \Delta x$ is obtained; the formula $dy = \Delta y$ is also obtained in the same way, therefore, $dz = \frac{\partial z}{\partial x}\Delta x + \frac{\partial z}{\partial y}\Delta y = \frac{\partial z}{\partial x}dx + \frac{\partial z}{\partial y}dy$ is obtained' [2,3]. Obviously it is considered that $z = x, z = y$ at the same time.

The definition of the equation $dx = \Delta x$ results in not only derivation of formula $dy = f(x_0)\Delta x = f(x_0)dx$, but also distortion of the general relationship between the change ($\Delta x$) and the differential ($dx$). The possibility exists in that $x$ is a function of another independent variable. In a function $y = F(x)$, $x$ is the reason and $y$ is the result. In turn, $y$ is also the reason of its sur-result and $x$ is also the result of its sub-reason. These relationships without beginning or end, existing in our world, are described as $z = E(y)$, $y = F(x)$, $x = G(t)$ algebraically. Therefore, $\Delta x \neq dx$ is obtained from the formula $\Delta x = g(t_0)\Delta t + o_g(\Delta t) = dx + o_g(\Delta t)$, except in special cases. Regarding a function of several variables, 'given the defined function $f(x) = f(x^1,\cdots,x^m)$ in the neighborhood of the point $x_0 = (x_0^1,\cdots,x_0^m)$, if the existence of $A_i$ ( $A_i \in R, i=1,\cdots,m$ ) leads to the equation $f(x_0+\Delta x) - f(x_0) = \sum_{i=1}^{m} A_i \Delta x^i + o(\|\Delta x\|)$ in a small neighborhood $\Delta x = (\Delta x^1,...,\Delta x^m)$, the function $f$ is differentiable and $\sum_{i=1}^{m} A_i \Delta x^i$ is called the total differential at the point $x_0$, expressed as $df(x_0) = \sum_{i=1}^{m} A_i \Delta x^i = \sum_{i=1}^{m} A_i dx^i$, where it is defined that $\Delta x^i = dx^i$, $i=1,......,m$' [4-6]. If each variable $x^i$ of the function $f(x) = f(x^1,...,x^m)$ is a differentiable function dependent on the variables $t = (t^1,...,t^k)$ (i.e., $x^i = x(t) = x^i(t^1,...,t^k)$ $i=1,......,m$), the equation $\frac{\partial}{\partial t^j}(f(x^1(t),...,f^m(t))) = \sum_{i=1}^{m} \frac{\partial f}{\partial x_i}(x(t)) \cdot \frac{\partial x_i}{\partial t_j}(t), j=1,...,k$ is obtained according to rules of derivation of composite functions, where $dx^i \neq \Delta x^i$.

Apparently whether $dx$ is equal to $\Delta x$ is related to the choice of a sign. Actually it is related to whether the differential is $\Delta y$ or the linear main-part of $\Delta y$. Otherwise, $o(\|\Delta x\|)$ could not be removed and it is impossible for the formation of Cauchy system. Thus it is not suitable to consider or define $dx = \Delta x$.

Even though the differential of the independent variable $x$ ($dx$) could be equal to $\Delta x$, cases regarding composite functions should be excepted. Has composite functions been dealt with in current calculus system? Sure! Differentials needs to be defined by derivatives and in turn, derivatives are defined by the differential quotient in current calculus system. The current definition of differentials leads to not only self-contradiction of the differentiation of composite functions, but also confusion in the derivation of composite functions, implicit functions, parametric equations and polar equations.

Indefinite integration could be achieved by integrating an objective expression ( $f(x)dx$ ) directly and indirectly (i.e., integration by parts), according to the Substitution Rule in most cases. Either the First Substitution Rule or the Second

Substitution Rule is deduced based on the concept of composite functions.

The similar problems also exist in the calculation of definite integrals, which should be calculated according to the corresponding indefinite integrals. The definite integration mentioned here, of course, is a concept, not the calculation of the sum in definite integrals in several cases.

If the mistakes mentioned above are made structurally, there exist more mistakes and defects in the meaning in the theories on current calculus as follows.

1. The definition of derivatives is not appropriate. A derivative should have been defined by differentials, but it is defined by $\lim_{\Delta x \to 0} \frac{\Delta y}{\Delta x}$ in a complicated way due to the wrong definition of differentials. The so-called derivative could be obtained by calculating the limitation of a rate of changes (i.e., the ratio of $\Delta y$ and $\Delta x$). However, the primitive function could not be obtained via a reverse way.

2. Although the part of indefinite integration is correct structurally, it is wrong logically in the meaning. A derived function $f(x)$ is restored as $F(x) = r_x f(x) + C$ if $r$ is defined as a sign of restoring and $x$ is an independent variable of restoring. Thus the formula $F'(x) = [r_x f(x)]' = f(x)$ is obtained, indicating that the derived function $f(x)$ is restored in the form of $r_x f(x)$ although the formula $d[r_x f(x)] = f(x)dx$ also exists. Similarly, the formula $d[\int f(x)dx] = f(x)dx$ indicates that the differentiation is restored in the form of $\int f(x)dx$ although the formula $[\int f(x)dx]' = f(x)$ also exists. The two completely different concepts of a family of primitive functions and difference of primitive functions are mixed up in the theories on current calculus. $r_x f(x)$ and $\int f(x)dx$ represent the family and difference of primitive functions although their derivatives and differentials are equal to each other, respectively. The difference of primitive functions is obviously misunderstood as the family of primitive functions if it is calculated by integration of $f(x)$ (i.e., $\int f(x)dx$).

3. The definite integral, defined by the limitation of sum of divided subareas, could be derived with the expression as $F(b) - F(a) = \int_a^b f(x)dx$ or $F(x) - F(a) = \int_a^x f(x)dx$, which are correct formally, but wrong logically in the meaning. If the differentiation is restored by summation, there is no need to divide.

Even it is necessary to divide, calculating the limitation does not make sense logically. Because a differential is a finite amount, a product by 'dividing'. The restoring is achieved just by summation and there is no need to calculating the limitation. The demonstration of Newton-Leibniz equation by Cauchy has been performed based on misunderstanding the difference of primitive functions as the family of primitive functions.

4. The principles of solutions according to the current definition of differentials are wrong for differential equations in the mathematical sense, not to mention differential equations containing composite functions in the structural sense. The differential $dx$, equal to $\Delta x$, is the change in the independent variable $x$; while the differential $dy$ is not the change ($\Delta y$) in the function $y$, but the linear main-part of $\Delta y$. Therefore, the value of $f(x)$ could be calculated by $\frac{dy}{\Delta x}$, but $\frac{dy}{\Delta x}$ is not consistent with $\lim_{\Delta x \to 0} \frac{\Delta y}{\Delta x}$ in the logic meaning. The formula $\Delta y = f(x)\Delta x + o(\Delta x)$ is deduced directly from the ensemble of $F(x)$; while $dy = f(x)\Delta x = f(x)dx$ is defined, not directly from the ensemble of $F(x)$. Thus it does not make sense by making use of the formula $dy = f(x)\Delta x = f(x)dx$ to solve differential equations in the meaning. In addition the difference and family of primitive functions are mixed up in the part of indefinite integrals.

There are more mistakes and defects in the theories on current calculus besides those mentioned above. However, either those mentioned above exists or one of them could be justified logically, to say the least, indicates that the theories on current calculus is not a 'rigorous and integrated system'.

## 2 Novel theories on calculus
### 2.1 Basic ideas

The definition of the real number and the point is suggested to be corrected as follows: the transformation of two equal real numbers to two different ones is subject to a Same-Different (S-D) transitional process; in the same way, the transformation of two different real numbers to two equal ones is subject to a Different-Same (D-S) transitional process; the transformation of two coincident points to two separated ones is subject to a Coincident-Separated (C-S) transitional process; in the same way, the transformation of two separated points to two coincident ones is subject to a Separated-Coincident (S-C) transitional process. In addition, it is defined that the extension at the corresponding level of an arbitrary real number is equal to zero; while the sub-extension (quasi-extension) is not equal to zero; the measurement at the corresponding level of a point is equal to zero; while the sub-measurement (quasi-measurement) is not equal to zero. The corrected definition of the real number and the point is considered as the novel quantity-figure model.

The essential difference between the traditional and novel quantity-figure models is the existence of the S-D (or D-S) transitional state of real numbers and C-S (or S-C) transitional state of points in the later model. Then the concept of true continuity is produced due to the existence of such transitional states. The assemblage, made up of numbers (or number groups) in the S-D (or D-S) transitional state between each element in a sequence, is defined as continuous variables. If the assemblage, produced by the mapping of continuous variables, is subject to the S-D (or D-S) transitional state between each element, such mapping relationship is defined as continuous functions.

Regarding a quantity $x_0$, whether it is an independent or dependent variable, there are three types of its increasing process $\Delta x = x - x_0$. The first, the increasing process is achieved by an accumulation of $\Delta x_i$ (a finite quantity) for n times (i.e., $\sum_{i=1}^{n} \Delta x_i = \Delta x$). The second, the process is achieved by an accumulation of $\delta x_i$ (infinitesimal) for $\infty$ times (i.e., $\sum_{i=1}^{\infty} \delta x_i = \Delta x$). The third, the interval, even less than $\delta x_i$, is equal to the 'difference' of two numbers (or points) in a S-D (C-S) transitional state. Thus the accumulated times are even larger than $\infty$, which is defined as super infinity expressed with the sign $\backsim$. If the interval is expressed as $dx_i$, the formulas $x_i - x_{i-1} = dx_i$ and $\sum_{i=1}^{\backsim} dx_i = \Delta x$ are obtained. In the three cases, the increment is actually divided into n, $\infty$ and $\backsim$ intervals, respectively. With the increase of the number of intervals, the measurement of each interval decrease from measurably finite quantity, unmeasurable infinitesimal quantity, to extreme quantity which is equal to zero from the level of $\Delta x$ due to the existence of D-S (or S-C) transitional state of the two values (or points). In such case $\delta x_i$ becomes so small that it is in the absent-present transitional state. Thus $dx_i$ is the 'bridge' for the transformation from absence to presence and unity of the opposite couple 'absence and presence'. This 'bridge' is called 'werden' in the book 'Logic' of Hegel. The interval $dx_i$ as the 'werden' is not only absence, but also presence; it is the absence in the presence and the presence in the absence. In a word, it is quasi-presence, which is consistent with the concept of 'werden'. The quantity of quasi-presence is equal to zero from the level of presence; while it is not equal to zero from the level of absence. The transitional state of numbers or points exists based on their characteristic of quasi-presence. For the traditional quantity-figure model, there are no S-D (or D-S) transitional states for numbers or C-S (or S-C) transitional states for points due to the lack of quasi-presence.

'Werden' can be translated into 'appearance' or 'disappearance', which can be united and evolve into another word 'change'. 'Werden' is the 'bridge' between absence and presence, which demonstrates specifically microcosmic relationship of analytic geometry. Figure 2.1 shows the geometrical schematic of an arc, a secant and a tangent. In Figure 2.1, as $\Delta x$ is becoming to be equal to $dx_i$, the arc $\overset{\frown}{AB}$, the secant $AB$ and the tangent $AC$ go into the state of 'werden', that is points A and B has entered the S-C transitional state. In this case the arc $\overset{\frown}{AB}$ and the secant $AB$ coincide with the tangent $AC$ in a way of reserving their (the arc and secant) restoring factors of relative quasi-presence (expressed with $o$). Thus all the quantities (or the drawings) enter the geometrical state without measurement, which is a state of topology, from the geometrical state with measurement. This accounts for that 'the replacement of the arc and the secant by the tangent' is accurate deduction, not approximate one. The corresponding microcosmic relationship is shown in Figure 2.2. Points A and B enter the S-C transitional state and point C coincides with B, with the corresponding algebraic relationship as follows.

$$dy \equiv F(x+dx) - F(x) \equiv \{tg[\theta(x)] + \alpha(dx)\}dx \equiv [f(x) + \alpha(dx)]dx \equiv f(x)dx +$$

$$o(dx) = f(x)dx \qquad (2.1)$$

Where the sign $\equiv$, which is a dialectic one and represents the relationship on mutual transformation and equality between two quantities having several levels (e.g., 'presence', 'quasi-presence'), is different from the sign =; $tg[\theta(x)]$ represents the slope function of $F(x)$ at the point $x_{i-1}$ and $\alpha(dx)$ is the relatively quasi-presence restoring function of $tg[\theta(x)]$; $o(dx) = \alpha(dx)dx$ is the relatively quasi-presence restoring factor of $F(x)$ at the point $x_{i-1}$, which is in the sub-level of $f(x)dx$.

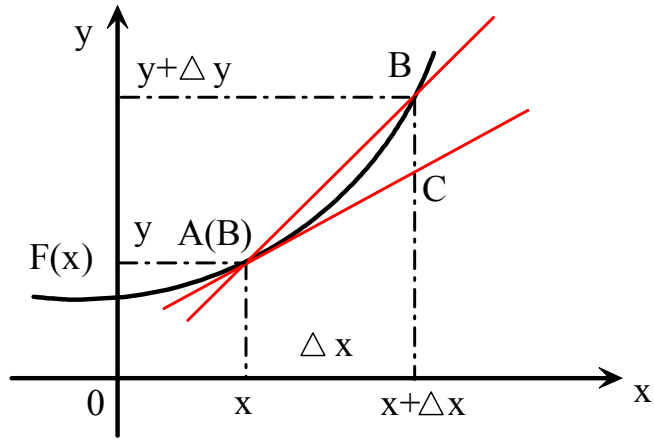

*Figure 2.1 Geometrical schematic of the arc $\overset{\frown}{AB}$, the secant $AB$ and the tangent $AC$.*

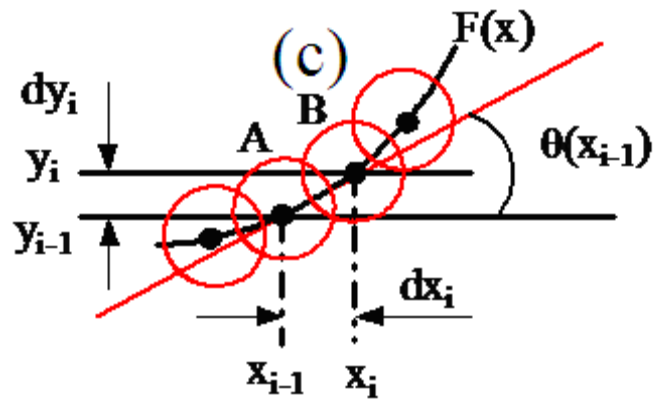

*Figure 2.2 Schematic of microcosmic relationship of the arc $\overset{\frown}{AB}$, the secant $AB$ and the tangent $AC$.*

Actually the theorem $\Delta y = f(x)\Delta x + o(\Delta x)$ exists based on the continuousness of $\dfrac{\Delta y}{\Delta x}$ in the interval $U(x,\eta)$. When $\Delta x$ becomes to be equal to $dx$, the formula $dy \equiv f(x)dx + o(dx) = f(x)dx$ is obtained.

**2.2 Definition of new concepts and the corresponding calculations**

The definition of differentials, derivatives, primitive functions and integrals and the corresponding calculations in the novel theories on calculus are introduced.

*2.2.1 Differentials and derivatives*

Given a function $y = F(x)$ (also applicable for a function of several variables), if the formula $dy \equiv [f(x)+\alpha(dx)]dx \equiv f(x)dx + o(dx) = f(x)dx$ exists in the interval $U(x,\eta)$, $dy \equiv [f(x)+\alpha(dx)]dx \equiv f(x)dx + o(dx)$ is defined as a dynamic differential function of $F(x)$ and $dy = f(x)dx$ as a static differential function of $F(x)$, called dynamic and static differentials in abbreviation, respectively.

At the same time, $F^o(x) \equiv \dfrac{dy}{dx} \equiv f(x) + \alpha(dx)$ is defined as a dynamic derived function of $F(x)$ and $F'(x) = \dfrac{dy}{dx} = f(x)$ as a static derived function of $F(x)$, called dynamic and static derivatives in abbreviation, respectively.

Two specific examples are given as follows.

EXAMPLE 1 Determine the dynamic and static differential functions as well as dynamic and static derived functions of $y = x^3$, respectively.

SOLUTION
Dynamic differential function
$$dy \equiv (x+dx)^3 - x^3 \equiv 3x^2 dx + 3x(dx)^2 + (dx)^3$$
$$\equiv [3x^2 + 3xdx + (dx)^2]dx$$

Static differential function
$$dy = 3x^2 dx$$

Dynamic derived function
$$F^o(x) \equiv \dfrac{dy}{dx} \equiv 3x^2 + 3xdx + (dx)^2$$

Static derived function
$$F'(x) = \dfrac{dy}{dx} = 3x^2$$

The structures of the dynamic differential and derived functions ($dy \equiv f(x)dx + o(dx)$ and $F^o(x) \equiv f(x) + \alpha(dx)$) have been clearly presented in this example with the obtained formulae $o(dx) \equiv 3x(dx)^2 + (dx)^3$ and $\alpha(dx) \equiv 3xdx + (dx)^2$.

EXAMPLE 2 Determine the second dynamic derived function of $y = x^3$.

SOLUTION $d_1 x$ and $d_2 x$ are used to represent the first and second differentials of x, respectively.

The first dynamic differential function
$$d_1 y \equiv 3x^2 d_1 x + 3x(d_1 x)^2 + (d_1 x)^3$$

The second dynamic differential function

$$d_2 d_1 y \equiv 3(x+d_2 x)^2 d_1 x + 3(x+d_2 x)(d_1 x)^2 + (d_1 x)^3 - 3x^2 d_1 x - 3x(d_1 x)^2 - (d_1 x)^3$$
$$\equiv 6x d_1 x d_2 x + 3(d_1 x)^2 d_2 x + 3 d_1 x (d_2 x)^2$$
$$\equiv [6x + 3 d_1 x + 3 d_2 x] \cdot d_1 x \cdot d_2 x$$

The second dynamic derived function

$$F^{00}(x) \equiv \frac{d_2 d_1 y}{d_1 x \cdot d_2 x} \equiv 6x + 3 d_1 x + 3 d_2 x$$

From the above example, it is concluded that the restoring function of the second dynamic derived function is a function of two variables. Generally speaking, the restoring function of n'th derived function is a function of n variables ($dx$).

2.2.2 *Primitive function*

If the formulae $F^o(x) \equiv f(x) + \alpha(dx)$ or $F'(x) = f(x)$ exists, $F(x) = r_x f(x) + C$ is called the family of primitive functions of $f(x)$ (primitive function in abbreviation). In the equation $r_x f(x)$ is called the base of the primitive function, $r$ is the sign of restoring, and $x$ is the variable of restoring. $r_x f(x)$ is the correction of traditional indefinite integral. $\int f(x) dx$ is considered as the conventional form of $r_x f(x)$ according to historical expression.

EXAMPLE 3 Determine the primitive function of $y = \sin^2 x$.

SOLUTION

$$r_x \sin^2 x = r_x \frac{1 - \cos 2x}{2} = r_x \frac{1}{2} - r_x \frac{1}{2} \cos 2x = \frac{1}{2} x - \frac{1}{4} r_{2x} \cos 2x = \frac{1}{2} x - \frac{1}{4} \sin 2x$$

$$\therefore F(x) = \frac{1}{2} x - \frac{1}{4} \sin 2x + C$$

2.2.3 *Integral*

The differential is 'change' (i.e., werden). The process of the generation of $\Delta x$ is $\sum_{i=1}^{\infty} dx = \Delta x$. But if $\Delta x$ is regarded as the origination and $dx_i \equiv \sum_{i=1}^{\infty} \Delta x$ is regarded as the generation of $dx_i$, $\sum_{i=1}^{\infty} dx_i$ is the restoration of the differential, which is defined as an integral.

Given $y = F(x)$,

$$F(x) - F(x_0) \equiv \Delta y \equiv \sum_{i=1}^{\infty} dy_i \equiv \sum_{i=1}^{\infty} [f(x_{i-1}) + \alpha(dx_i)] dx_i$$
$$\equiv \sum_{i=1}^{\infty} [f(x_{i-1}) dx_i + o(dx_i)] = \sum_{i=1}^{\infty} f(x_{i-1}) dx_i \quad (2.2)$$

Where a sign of continuous summation ($\int_{x_0}^{x}$) is introduced as a equivalence of the sign $\sum_{i=1}^{\infty}$, which is called a symbol of pseudo-discrete summation. The equation (2.2) is converted to the Newton-Leibniz equation as follows.

$$F(x) - F(x_0) \equiv \int_{x_0}^{x} [f(x) + \alpha(dx)] dx = \int_{x_0}^{x} f(x) dx \quad (2.3)$$

Where $\int_{x_0}^{x}$ is summation of $[f(x) + \alpha(dx)] dx$ or $f(x) dx$ from $x_0$ to $x$ with an interval $dx$. $f(x)$ and $\alpha(dx)$ synchronize with $x$ which starts from $x_0$ throughout the process of summation.

The result of calculating the dynamic and static differentials of the Newton-Leibniz equation is

$$D\left\{\int_{x_0}^{x} [f(x) + \alpha(dx)] dx\right\} \equiv D[F(x) - F(x_0)] \equiv [f(x) + \alpha(dx)] dx$$

$$d\left\{\int_{x_0}^{x} [f(x) + \alpha(dx)] dx\right\} = d[F(x) - F(x_0)] = f(x) dx$$

Thus it is clear that the integral is the restoration of the dynamic and static differential.

Suppose $F(x)$ and its static derivative $f(x)$ is continuous on the closed interval $[x_0, x]$ where $x_0, x \in [a, b]$. Given a finite quantity $\ell$, $\ell$ points are chosen from the interval $(x_0, x)$ arbitrarily. $\infty' = \infty - \ell$ is defined and the following equation exists.

$$\int_{x_0}^{x}[f(x)+\alpha(dx)]dx \equiv \sum_{j=1}^{\infty}[f(x_{i-1})+\alpha(dx_i)]dx_i$$

$$\equiv \sum_{j=1}^{\infty'}[f(x_{j-1})+\alpha(dx_j)]dx_j + \sum_{k=1}^{\ell}[f(x_{k-1})+\alpha(dx_k)]dx_k$$

$$\equiv \int_{x_0}^{x}{}_l[f(x)+\alpha(dx)]dx$$

$$= \int_{x_0}^{x}{}_l f(x)dx$$

Where $\sum_{j=1}^{\infty'}[f(x_{j-1})+\alpha(dx_j)]dx_j \equiv \int_{x_0}^{x}{}_l[f(x)+\alpha(dx)]dx$ is equal to $\sum_{j=1}^{\infty}[f(x_{i-1})+\alpha(dx_i)]dx_i \equiv \int_{x_0}^{x}{}_l[f(x)+\alpha(dx)]dx$ in the quantity. But the former is lack of $\sum_{k=1}^{\ell}[f(x_{k-1})+\alpha(dx_k)]dx_k$ that represents a trend of change in microcosmic scale compared with the latter. For the limited growth of $\sum_{k=1}^{\ell}$, $\sum_{k=1}^{\ell}[f(x_{k-1})+\alpha(dx_k)]dx_k$ is equal to zero in the level of $f(x)$. The geometrical significance of $\sum_{k=1}^{\ell}[f(x_{k-1})+\alpha(dx_k)]dx_k$ is the generating trend of $\ell$ pieces of line manifold and $\ell$ pieces of surface manifold. The $\int_{x_0}^{x}{}_l[f(x)+\alpha(dx)]dx$ in lack of $\sum_{k=1}^{\ell}[f(x_{k-1})+\alpha(dx_k)]dx_k$ is called hypo-generation-integral.

If f(x) has $\ell$ points of the first discontinuity in the interval $(x_0, x)$ and F(x) is continues or has removable discontinuity in the interval $(x_0, x)$, then

$$\int_{x_0}^{x}f(x)dx = \int_{x_0}^{x}[f(x)+\alpha(dx)]dx = \int_{x_0}^{x}{}_\ell[f(x)+\alpha(dx)]dx$$
$$= \int_{x_0}^{x}{}_\ell f(x)dx = F(x)-F(x_0) \quad (2.4)$$

## 3. Revelations from the History of Calculus

Throughout the history of calculus, it is found that many brilliant thoughts of Newton and Leibniz have been misjudged; the Newton-Leibniz theories on calculus have returned in the form of the novel theories demonstrated in this paper; the current theories on calculus have deviated from the development of calculus since its establishment by Cauchy.

Newton and Leibniz invented calculus in 1667 independently, respectively. Their theories on calculus were developed and enriched by later mathematicians, with Euler being the representative [7].

Newton and Leibniz substituted an arc or a secant by a tangent and used the concept of a characteristic triangle, following the previous research. Newton said, 'The final ratio of any two quantities of the arc, secant and tangent is equal to each other. [8]' Leibniz said, 'A characteristic triangle is one that it retains the form of a triangle after the all the corresponding quantities have been removed. [9]' Of course, these are based on the deduction of the theory. Leibniz also did contributions in the establishment of the theory. He presented the formulae $dx^e = ex^{e-1}dx$ and $\int x^e dx = \frac{x^{e+1}}{e+1}$ in 1667 as well as $\int_a^b ydx = Z(b) - Z(a)$ in 1677. He gave the following explanation concerning the symbols, '…dx represents the difference of two adjacent x…A differential is like the contact angle of Euclid, which is smaller than any given quantity, but not zero…it is a qualitative zero. [8]' He also showed that 'd' is the reverse of '$\int$' in the meaning, and that '$\int$' is the extension of the initial letter of 'sum' [9]. Actually Leibniz had begun to try to reestablish the quantity-figure model. He said, 'For a point, the extension of it is equal to zero, not its measurement.' Indeed, the ideas of Newton and Leibniz are not consistent, and the correction by Leibniz on the definition of the point may not be correct. However, no one could deny that the ideas of Newton and Leibniz are embodied in the above mentioned description.

Bernard Bolzano and Cauchy criticized that the ideas of Newton and Leibniz were unstable, and that their theories were not clear, but this definitely is not responsible for the deviation of calculus from its right development. In fact, Cauchy could not understand the brilliant ideas of Newton and Leibniz unless he got rid of the thinking pattern of Jean d'Alembert. The 'zero' of Newton (in some cases it is zero; in some cases it is not [10]) and the 'qualitative zero' of Leibniz are embryos of 'werden'. However, d'Alembert criticized, 'A quantity could just be either in a state of presence or absence. If it exists, it has not disappeared, and vice versa. If a transitional state exists, it corresponds to a chimae consisting of the head of a lion, the body of a sheep and the tail of a sneak. [8]' Under such thinking pattern, it is no wonder that the theories on calculus established by Newton and Leibniz cannot be understood.

The application of the limit or the inequation is necessary for establishing the theories on calculus once the differential is chosen as a finite quantity. Irrational means are employed if problems cannot be treated by applying the limit and the inequation, such as the transformation of f(x) Δx into f(x)dx. However, the derivative f(x) cannot be clearly explained even if f(x)dx exists, which causes that the definition of the derivative precedes that of the differential. Consequently, the definition of the derivative by calculating the limitation of a rate of changes (i.e., the ratio of $\Delta y$ and $\Delta x$) is applied rather than by the inequality due to the simplicity of the former one. However, such definitions of the derivative and the differential result in the loss of the restoring function $\alpha(dx)$ and the restoring factor o(dx), and the impossibility of transformation from derivatives to primitive functions reversely. Actually the differential can be easily defined by the derivative, as is shown in the equation $dy \equiv \frac{dy}{dx} \cdot dx \equiv [f(x) + \alpha(dx)]dx = f(x)dx$. Cauchy avoided the concept of inverse differentials when he defined the definite integral, but he returned to primitive functions after he "demonstrated" the equation of Newton-Leibniz. As long as Cauchy admitted that the definite integral was the inverse differential, the logic defects in the definition of definite integrals by calculating the limitation of sum of divided subareas would be exposed. The equation of Newton-Leibniz would not be obtained if Cauchy did not employ the primitive function.

The equation of Newton-Leibniz could have been easily proved according to $F(b) - F(a) = \int_{F(a)}^{F(b)} dy = \int_{F(a)}^{F(b)} \frac{dy}{dx} \cdot dx = \int_a^b f(x)dx$. In current system, the equation of Newton-Leibniz could be obtained by the differential as a finite quantity, which is achieved via the following steps. Given $\Delta y_i = f(x_i)\Delta x_i + o(\Delta x_i)$, $F(b) - F(a) = \sum_{i=0}^{n-1} f(x_i)\Delta x_i + \sum_{i=0}^{n-1} o(\Delta x_i)$ can be obtained. Calculating the limits of both sides, the equation $F(b) - F(a) = \lim_{n \to \infty} \sum_{i=0}^{n-1} f(x_i)\Delta x_i + \lim_{\Delta x \to 0} \sum_{i=0}^{n-1} o(\Delta x_i) = \lim_{\Delta x \to 0} \sum_{i=0}^{n-1} f(x_i)\Delta x_i$ is obtained, where, of course, it is difficult to show $\lim_{\Delta x \to 0} \sum_{i=0}^{n-1} o(\Delta x_i) = 0$ directly.

The theories on current calculus presented by Cauchy, Riemann, Weierstrass, Darboux and Lebesgue, have a negative effect in practice. Before Cauchy, the concepts of Δx, Δy, dx and dy are independent, different from Δx and Δy limited by $\lim_{\Delta x \to 0} \frac{\Delta y}{\Delta x}$; different from dx and dy which are apart from directly deduced ensemble of F(x) due to the existence of dy=f(x)Δx=f(x)dx. Therefore, the Lagrange equation can be deduced reasonably. Siméon Poisson, making significant achievements in several fields, disagreed with Cauchy. He emphasized that the differential should be considered as an infinitesimal quantity, not a finite quantity. He said, 'the infinitesimal quantity…smaller than any given quantity with the identical quality. [8]' However, it is not equal to zero and attentions should be paid to 'the identical quality'. Actually the mentioned infinitesimal quantity is the 'Werden' corresponding to 'quasi-presence'.

In 1960s, Abraham Robinson 'proved that the structure of real numbers could be expanded to include infinitesimals and infinity by the model theory… providing logically reliable basis for the controversial concept of infinitesimals proposed by Leibniz logically for the first time [11]'. Robinson said that, 'The book (Non-standard Analysis) justifies the ideas of Leibniz.'   Therefore, the differential of 'a qualitative zero' is justified; while the self-contradictory differential of 'a finite quantity' by Cauchy and the corresponding theories on calculus are not justified.

The respectable mathematician Kurt Gödel in the field of mathematics and philosophy considered that the theories on calculus of Cauchy-Lebesgue system would end certainly. He said, 'The non-standard analysis in various forms will become the analytics in the future. [12]'

**4 About mathematics**

In this paper we have always managed to combine the theory with practice, better to say to demonstrate the activity of the novel theories on calculus than to justify them, when presenting the novel quantity-figure model and the novel theories. Because, in general, mathematics is formal science and should be subject to logical self-consistency.

The current theories on calculus have been challenged due to their logical self-contradiction in the structure and meaning; we are convinced of the novel quantity-figure model and theories on calculus for their self-justification.

Mathematics, similar to other subjects, refers to the approximation of human's brains to the 'truth' making use of models or constants in a way of logical system. In such case the mentioned 'truth' is about the quantity, figure and self-structure. Since the approximation is carried out based on models, the approximate points and degrees are dependent on the establishment of models, which is the most complicated problem during the constructing of scientific models although some of them appears to be simple. The constructer is required to be equipped with the ability to perceive, at least, and to predict the trend of the corresponding field from a strategical standpoint. Thus a remarkably representative and concise model could be established. However, the established model, no matter how representative it is, could only summarize some aspects of the real world or basic properties of substances, not all the aspects. Therefore, the characteristic of the established model should be strengthened, no matter how successful it is, and the challenges based on specific representativeness are rejected after its establishment. It is better to emphasize the sub-principles from the standpoint of axioms than theorems during the establishment of models.

Regarding the axioms, they should be subject to the principles of the novel formal logic although they could possibly not be deduced from the novel formal logic, which refers to the corrected formal logic. Actually, the antinomies occur due to the lack of laws of time limit and the transitional state during the transformation of substances could not be appropriately dealt with in the current system of formal logic. H. Poincaré said, 'Mathematics is generated from experience, but it should not be discriminated by the factors irrelevant to the experience. [8]'

The common point of mathematics and natural science is that both of them are approximate system to the 'truth'; while the difference is that mathematics is purely formal science. Therefore, a system could become the most successful mathematical science under the conditions including logical self-consistency, wide-ranging application and relatively best approximation to laws of quantities, figures and self-structure.

**5 Conclusions**

The basis for logical self-contradiction of the theories on the current calculus is exposed. A new quantity-figure model as the premise of mathematics is formed after the correction of definition of the real number and the point. A novel calculus system is established based on the concept of 'werden' and the theories on calculus are reestablished.